\documentclass[12pt,a4paper]{article}
\usepackage{amsmath,amsfonts,amssymb,amsthm,amscd}
\newtheorem{thm}{Theorem}
\newtheorem{lem}[thm]{Lemma}
\newtheorem{prop}[thm]{Proposition}

\newtheorem{cor}[thm]{Corollary}
\newtheorem{rem}[thm]{Remark}

\begin{document}

\title{A characterization of quaternionic projective space by
the conformal-Killing equation}

\author{Liana David and Massimiliano Pontecorvo}

\maketitle

{\bf Abstract:} We prove that a compact quaternionic-K\"{a}hler
manifold of dimension $4n\geq 8$ admitting a conformal-Killing
$2$-form which is not Killing, is isomorphic to the quaternionic
projective space, with
its standard quaternionic-K\"{a}hler structure.\\

\section{Introduction}

The existence of a non-constant smooth function on a Riemannian manifold,
which satisfies a well-chosen
differential equation can, under certain conditions, determine
the Riemannian manifold. Such characterizations in terms
of solutions of differential equations exist for
the sphere,
the complex projective space and the quaternionic projective space
(with their standard metrics), see \cite{blair}, \cite{obata},
\cite{ika}.
In this note we develop an alternative characterization
of the
quaternionic projective space using the conformal-Killing equation.

Recall that a $p$-form $u$ on a Riemannian manifold
$(M^{m}, g)$ is Killing, if its covariant derivative $\nabla u $  with
respect to the
Levi-Civita connection $\nabla$ is totally skew-symmetric, or,
equivalently, if $\nabla u =\frac{1}{p+1}du.$
Killing forms are natural generalizations of Killing vector fields
and were introduced for the first time by K. Yano in \cite{yano}. More
generally, one can consider conformal-Killing forms,
which are forms $\psi\in \Omega^{p}(M)$ satisfying the conformal-Killing
equation:
\begin{equation}\label{confkil}
\nabla_{Y} \psi  =\frac{1}{p+1} i_{Y}d\psi  -\frac{1}{m-p+1}Y\land
\delta \psi ,\quad\forall Y\in TM.
\end{equation}
(Here and everywhere in this note we identify tangent vectors with
their dual $1$-forms, using the Riemannian metric). Note that a
conformal-Killing form is Killing if it is coclosed. There is an
intensive literature on conformal-Killing forms. On K\"{a}hler
manifolds conformal-Killing forms are closely related to
Hamiltonian $2$-forms and have been completely classified (in the
compact case) in \cite{kahler}. In particular, conformal-Killing
forms exist on Bochner-flat K\"{a}hler manifolds and on
conformally Einstein K\"{a}hler manifolds. In this paper, we study
conformal-Killing $2$-forms on quaternionic-K\"{a}hler manifolds,
which are Riemannian manifolds with holonomy group included in
$\mathrm{Sp}(n)\mathrm{Sp}(1)$. It is known that on a compact
quaternionic-K\"{a}hler manifold any Killing $p$-form (with $p\geq
2$) is parallel \cite{qK}. For this reason, we shall be interested
in conformal-Killing $2$-forms which are not Killing. Our aim is
to prove the following result:

\begin{thm}\label{main}\begin{enumerate}\item A compact,
connected, quaternionic-K\"{a}hler manifold $(M, g)$ of dimension
$4n\geq 8$ admits a conformal Killing $2$-form which is not
Killing if and only if it is isomorphic to the quaternionic
projective space $\mathbb{H}P^{n}$, with its standard
quaternionic-K\"{a}hler structure.

\item Let $g_{\mathrm{can}}(\nu )$ be the standard metric of
$\mathbb{H}P^{n}$, with reduced scalar curvature $\nu >0.$ The map
which associates to a Killing vector field $X$ on
$(\mathbb{H}P^{n},g_{\mathrm{can}}(\nu ))$ the $2$-form
\begin{equation}\label{forma}
\psi  :=-\frac{2}{\nu (4n-1)}(\nabla X)^{S^{2}H}
+\frac{4}{\nu (4n-1)}(\nabla X)^{S^{2}E}
\end{equation}
is an isomorphism from the space of Killing vector fields to the
space of conformal-Killing $2$-forms on
$(\mathbb{H}P^{n},g_{\mathrm{can}}(\nu ))$, with inverse the
codifferential: $\delta (\psi ) =X.$

\end{enumerate}

\end{thm}

The plan of the paper is the following. In Section \ref{prel} we
recall some basic facts about quaternionic-K\"{a}hler manifolds
and the conformal-Killing equation. An important feature for us is
that the codifferential of a conformal-Killing $2$-form on an
Einstein manifold (hence also on a quaternionic-K\"{a}hler
manifold) is a Killing vector field \cite{semc}. In Section
\ref{comp} we determine the general form of conformal-Killing
$2$-forms on a compact quaternionic-K\"{a}hler manifold $(M, g)$
of dimension $4n\geq 8$ (see Proposition \ref{wx}). We show that
there are no conformal-Killing, non-Killing $2$-forms on $(M, g)$,
unless the reduced scalar curvature $\nu$ is positive, in which
case a Killing vector field $X$ on $(M, g)$ is the codifferential
of a conformal-Killing $2$-form if and only if it belongs to the
kernel of the quaternionic-Weyl tensor $W$ of $(M, g).$ In Section
\ref{kerw} we conclude the proof of Theorem \ref{main}, by showing
that if $(M, g)$ (still compact, quaternionic-K\"{a}hler, with
$\nu >0$ and of dimension $4n\geq 8$) admits a (non-trivial)
Killing vector field $X$ in the kernel of $W$, then it is
isometric to the standard quaternionic projective space (see
Proposition \ref{kil}). Our method to prove this statement is to
consider the Hamiltonian function $f^{X}$ of the natural lift
$X^{Z}$ of $X$ to the twistor space $Z$ of $(M, g)$ and to show
that it satisfies a certain differential equation introduced by
Obata in \cite{obata}. By a result of Obata (see \cite{obata},
Theorem $C$) this implies that $Z$ (with its standard
K\"{a}hler-Einstein structure) is isomorphic to the complex
projective space, with its Fubini-Study metric and therefore the
quaternionic-K\"{a}hler manifold $(M, g)$ must be isomorphic to
$(\mathbb{H}P^{n},g_{\mathrm{can}}(\nu ))$.

Similar type of results appear in the literature. In four
dimensions, a conformal oriented manifold with degenerate and
coclosed self-dual Weyl tensor is necessarily anti-self-dual (see
\cite{der}, page 454). In the same framework, gradient
quaternionic vector fields on quaternionic-K\"{a}hler manifolds
lie in the kernel of the quaternionic-Weyl tensor; moreover, if
the quaternionic-K\"{a}hler manifold is compact and non
Ricci-flat, then it has no non-zero gradient vector fields, unless
it is isomorphic to the quaternionic projective space, with its
standard quaternionic-K\"{a}hler structure, see \cite{gradient},
\cite{gradient1}.

In the last Section of the paper we determine the dimension of the
space of conformal-Killing $2$-forms defined on a compact,
quaternionic-K\"{a}hler manifold. This is a
consequence of Proposition \ref{wx} of Section \ref{comp}.\\

{\bf Acknowledgements:} This research has been undertaken while one
of the authors, Liana David, was a junior visiting research fellow
at Centro di Ricerca Matematica "Ennio de Giorgi", Scuola Normale
Superiore di Pisa, Italy. Both authors thank Centro di Ricerca for
financial support and excellent working conditions. We also thank
Vestislav Apostolov, Paul Gauduchon, Andrei Moroianu, Simon Salamon
and Uwe Semmelmann for their interest in this work.

\section{Preliminary material}\label{prel}

In this paper we shall use the following conventions and
notations. All our manifolds will be smooth and connected. The
space of smooth sections of a vector bundle (real or complex)
$V\rightarrow M$ over a manifold $M$ will be denoted by $\Gamma
(V)$. As in Introduction, $\Omega^{k}(M)$ will denote the space of
smooth, real-valued $k$-forms on $M$. Finally, all our
quaternionic-K\"{a}hler manifolds will be of (real) dimension bigger or equal to eight.\\

{\bf Quaternionic-K\"{a}hler manifolds.} Let $(M, g)$ be a
quaternionic-K\"{a}hler manifold, i.e. a Riemannian manifold with
holonomy group included in $\mathrm{Sp}(n)\mathrm{Sp}(1).$ The
endomorphism bundle of $TM$ has a distinguished, parallel (i.e.
preserved by the Levi-Civita connection) rank $3$-subbundle $Q$,
called the quaternionic bundle, which is locally generated by a
system of three almost complex structures $\{ J_{1}, J_{2},
J_{3}\}$ (called an admissible basis of $Q$) subject to the
quaternionic relations
$$
J_{1}^{2} = J_{2}^{2} = J_{3}^{2} =- \mathrm{Id},\quad J_{i}
J_{j}= - J_{j}J_{i},\quad\forall i\neq j.
$$

Like for conformal $4$-manifolds, there are two locally defined
complex vector bundles $H$ and $E$ over $M$, of rank $2$ and $2n$
respectively, associated to the standard representations of
$\mathrm{Sp}(1)$ and $\mathrm{Sp}(n)$ on $\mathbb{C}^{2}$ and
$\mathbb{C}^{2n}$. The bundles $E$ and $H$ play the role of the
spin bundles in conformal geometry. In particular, the
complexification $T_{\mathbb{C}}M$ is isomorphic to the tensor
product $E\otimes H$ and the complexification of the bundle of
$2$-forms has the parallel decomposition
\begin{equation}\label{decf}
\Lambda^{2}(T^{*}_{\mathbb{C}}M) = S^{2}H\oplus S^{2}E\oplus
(S^{2}H\otimes  \Lambda_{0}^{2}E),
\end{equation}
where $\Lambda_{0}^{2}E\subset \Lambda^{2}E$ is the kernel of the natural
contraction with the
standard symplectic form on $E$.
The bundle $S^{2}H$ is isomorphic to the complexification of
the bundle $Q$, and
$S^{2}E$ is isomorphic to the complexification of the bundle of
$Q$-hermitian forms, i.e. $2$-forms $\psi$ which satisfy
$$
\psi (AX, Y) = - \psi (X, AY),\quad\forall A\in Q, \quad \forall X, Y\in TM.
$$
For a $2$-form $\psi \in\Omega^{2}(M)$, we denote by
$\psi^{S^{2}H}$, $\psi^{S^{2}E}$ and $\psi^{S^{2}H\otimes\Lambda^{2}_{0}E}$
its projections on the three factors of the decomposition (\ref{decf}).
Note that, when $\psi = X\land Y$ is decomposable,
$$
(X\land Y)^{S^{2}H} =\frac{1}{2n}\sum_{i=1}^{3}\omega_{i}(X, Y)\omega_{i}
$$
and
$$
(X\land Y)^{S^{2}{E}}= \frac{1}{4}\left( X\land Y +\sum_{i=1}^{3}J_{i}X\land J_{i}Y\right) ,
$$
with respect to any admissible basis
$\{ J_{1}, J_{2}, J_{3}\}$
of $Q$, with associated K\"{a}hler forms $\omega_{i}= g(J_{i}\cdot ,\cdot ).$

We now turn to the curvature of the quaternionic-K\"{a}hler manifold $(M, g).$
The metric $g$ is Einstein and its curvature has the expression
\begin{equation}\label{rg}
R^{g}(X, Y) =-\frac{\nu}{4}\left( X\land Y +\sum_{i=1}^{3}J_{i}X\land J_{i}Y
+2\sum_{i=1}^{3}\omega_{i}(X, Y)\omega_{i}\right) +W(X, Y)
\end{equation}
where $\nu :=\frac{k}{4n(n+2)}$ is the reduced scalar curvature
($k$ being the usual scalar curvature) and $W$ is the
quaternionic-Weyl tensor, which is a symmetric endomorphism of
$S^{2}E$ and is in the kernel of the Ricci contraction. The tensor
$W$ plays the role of the Weyl tensor in conformal geometry, in
the sense that if $W=0$, then $(M,Q)$ is locally isomorphic, as a
quaternionic manifold,
to the quaternionic projective space $\mathbb{H}P^{n}.$\\

There is one more piece of information we need to recall, namely the
twistor space of $(M, g).$
The bundle $Q$ has a natural Euclidian metric
$\langle \cdot ,\cdot \rangle$, for which any admissible
basis is orthonormal.
The twistor space of $(M, g)$
is the total space $Z$ of the unit sphere bundle
of $Q$, i.e. the set of all complex structures of tangent spaces
of $M$, which, seen as endomorphisms of $TM$, belong to $Q$.
The fibers $Z_{p}:= \pi^{-1}(p)$
of the twistor projection $\pi :Z\to M$, called
twistor lines, are complex manifolds, with complex structure
$\mathcal J$ defined by
\begin{equation}\label{cs}
{\mathcal J}(A):= J\circ A,\quad\forall A\in T_{J}Z_{p},\quad\forall J\in Z_{p},\quad\forall p\in M.
\end{equation}
Note that $\mathcal J$ is well-defined, since
$$
T_{J}Z_{p}= \{ A\in Q_{p}:\quad A\circ J+J\circ A =0\} = J^{\perp}\subset Q_{p},
$$
where $\perp$ denotes the orthogonal complement with respect to the
metric $\langle \cdot ,\cdot\rangle .$ The twistor space $Z$
has a standard (integrable) complex structure, also denoted by $\mathcal J$,
and, when
$\nu >0$, a K\"{a}hler-Einstein metric $\bar{g}$.
In order to define $\mathcal J$ and $\bar{g}$,
consider the horisontal bundle $H^{\nabla}\subset TZ$
associated to the Levi-Civita connection $\nabla$,
acting on the twistor bundle $\pi :Z\to M$.
The complex structure $\mathcal J$ preserves $H^{\nabla}$ and
the twistor lines. Its restriction to $H^{\nabla}_{J}$
(for any $J\in Z$) is defined tautologically, using the
linear isomorphism $\pi_{*}:H^{\nabla}_{J}\to T_{p}M$ (where
$p:= \pi (J)$)
and its restriction to the twistor lines coincides with the standard
complex structure of the twistor lines, defined in
(\ref{cs}). The metric $\bar{g}$ is defined in the following way:
on $H^{\nabla}$ it is the pull-back of $g$; the twistor lines are
$\bar{g}$-orthogonal to $H^{\nabla}$;
when $\nu =1$
the restriction of $\bar{g}$ to the twistor lines is the standard
metric on $S^{2}$ of curvature one; equivalently,
the restriction of $\bar{g}$ to a twistor line $Z_{p}$
is induced by the Euclidian metric $\langle \cdot ,\cdot\rangle $ of the fiber
$Q_{p}$ of $Q$ over $p$, by means
of the inclusion $Z_{p}\subset Q_{p}$.
The twistor projection $\pi :(Z,\bar{g})\to (M,g)$ becomes a Riemannian submersion
with totally geodesic fibers.\\

{\bf The conformal-Killing equation.} The conformal-Killing
equation on $p$-forms on a compact Riemannian manifold $(M^{m}, g)$
can be written in the alternative form \cite{semc}
\begin{equation}\label{inplus}
q(R)\psi = \frac{p}{p+1}\delta d\psi +\frac{m-p}{m-p+1}d\delta \psi ,
\end{equation}
where $q(R)$ is a bundle endomorphism of $\Lambda^{p}(T^{*}M)$,
related to the Laplacian
$\Delta = d\delta +\delta d$ by  the formula
$\Delta = \nabla^{*}\nabla  +q(R)$, where, for any form $\psi$,
$\nabla^{*}\nabla\psi  = -\sum_{i} \nabla^{2}(\psi )(E_{i}, E_{i})$,
$\{ E_{i}\}$ is an orthonormal local frame of $TM$
and $\nabla^{2}(\psi )(X, Y) :=\nabla_{X}\nabla_{Y}\psi  -  \nabla_{\nabla_{X}Y}\psi $, for any
vector fields
$X$ and $Y$. More explicitly,
\begin{equation}\label{q}
q(R)(\psi ):= \sum_{i,j=1}^{m} E_{j}\land i_{E_{i}}R^{g}(E_{i},
E_{j})(\psi ),\quad\forall \psi \in \Lambda^{p}(T^{*}M).
\end{equation}
In (\ref{q}) $R^{g}$ is the curvature operator of the Levi-Civita
connection acting on the form bundle, so that $R^{g}(E_{i},
E_{j})$ is an endomorphism of $\Lambda^{p}(T^{*}M)$ and
$R^{g}(E_{i}, E_{j})(\psi)$ denotes the action of $R^{g}(E_{i},
E_{j})$ on the $p$-form $\psi $. An important feature of the
curvature operator $q(R)$ is that it preserves the parallel
subbundles of the form bundle (see, for example, \cite{kahler}).
If $(M, g)$ is a symmetric space, then the operator $q(R)$ is
parallel and commutes with the Laplace operator $\Delta$. If,
moreover, $M$ is compact, $\Delta$ acts on the (finite
dimensional, see \cite{semc}) vector space of conformal-Killing
forms on $(M, g)$ and is diagonalisable on this space. In
particular, any conformal-Killing form defined on a compact
symmetric space can be written as a linear combination of
conformal-Killing forms, which are also eigenforms of the Laplace
operator. (It is expected that this observation, together with the
estimates found in \cite{semqk} on the eigenvalues of Laplace
operator on a compact quaternionic-K\"{a}hler manifold with
positive scalar curvature, might be useful to understand higher
degree conformal-Killing forms on Wolf spaces; further
investigation in this direction
is needed).\\

Suppose now that
$(M^{4n}, g)$ is a
quaternionic-K\"{a}hler manifold, with reduced scalar curvature
$\nu .$ We shall be interested in the conformal-Killing equation on
$2$-forms defined on $(M, g)$.
The operator $q(R)$ acts on $S^{2}H$ and $S^{2}H\otimes
\Lambda_{0}^{2}E$ by scalar multiplication
as follows (see \cite{semqk}, Lemma 2.5):
\begin{equation}\label{constante}
q(R)\vert_{S^{2}H} = 4\nu\mathrm{Id},\quad
q(R)\vert_{S^{2}H\otimes \Lambda^{2}_{0}E} = 2\nu
(n+2)\mathrm{Id}.
\end{equation}
(We remark that the operator $q(R)$ of the conformal-Killing
equation (\ref{inplus}) differs by a multiplicative factor of $2$
from the operator $q(R)$ considered in Lemma 2.5
of \cite{semqk}). The action of $q(R)$ on $2$-forms
still preserves $S^{2}E$, but it is not in general a scalar action
on this space.

The metric $g$ being Einstein, the codifferential of a
conformal-Killing $2$-form on $(M, g)$ is a Killing vector field
(this is an easy consequence of Proposition 5.2 of \cite{semc}).

\section{Conformal Killing $2$-forms on
quaternionic-K\"{a}hler manifolds}\label{comp}

On a quaternionic-K\"{a}hler manifold with non-zero scalar
curvature, the codifferential $\delta$ is a linear isomorphism
from the space of sections of $S^{2}H$ which are solutions of the
twistor equation to the space of Killing vector fields (see
\cite{inventiones}, Lemma 6.5 and \cite{amp1}, Proposition 5.6).

In this Section we will prove an analogous statement for
conformal-Killing $2$-forms. As already mentioned in the
Introduction and in Section \ref{prel}, the codifferential
$\delta$ sends conformal-Killing $2$-forms on a
quaternionic-K\"{a}hler manifold $(M, g)$ to Killing vector fields
and its kernel is the space of Killing $2$-forms. We will now
determine the image of $\delta$ (considered as a map from
conformal-Killing $2$-forms), under the additional assumption that
$M$ is compact. More precisely, we prove the following result.

\begin{prop}\label{wx} Let $\psi$ be a conformal-Killing $2$-form
on a compact quaternionic-K\"{a}hler manifold $(M, g)$ of
dimension $4n\geq 8$ and reduced scalar curvature $\nu$. If $\nu
>0 $ then
\begin{equation}\label{f}
\psi =-\frac{2}{\nu (4n-1)}(\nabla X)^{S^{2}H}
+\frac{4}{\nu (4n-1)}(\nabla X)^{S^{2}E}
+ u,
\end{equation}
where $X:=\delta (\psi)$ is the codifferential of $\psi$ and
$u\in\Omega^{2}(M)$ is parallel. Moreover, $W(X, \cdot )=0$, where
$W$ denotes the quaternionic-Weyl tensor of $(M, g).$ If $\nu \leq
0$ then $\psi$ is parallel.
\end{prop}

\begin{rem}\label{remi}{\rm For Killing forms, Proposition \ref{wx} reduces to the statement
proved by Moroianu and Semmelmann in \cite{semqk}, namely that any
Killing $2$-form on a compact quaternionic-K\"{a}hler manifold is
parallel. For this reason, Proposition \ref{wx} is relevant when
$\psi$ is a conformal-Killing, but not Killing $2$-form.}
\end{rem}

We now prove Proposition \ref{wx}. The case $\nu \leq 0$ is an
easy consequence of the following observations: there are no
(non-trivial) Killing vector fields on a compact
quaternionic-K\"{a}hler manifold $(M, g)$, with $\nu <0$ (this is
an application of the Weinzenb\"{o}ck formula, see also
\cite{der}, Theorem 1.84). Similarly, if $\nu =0$ then any Killing
vector field on $(M, g)$ is harmonic, and, if coexact, it is
identically zero ($M$ being compact). Due to these facts, any
conformal-Killing $2$-form on a compact, quaternionic-K\"{a}hler
manifold  with non-negative scalar curvature is Killing, hence
parallel \cite{semqk}. This proves Proposition \ref{wx} when $\nu \leq 0.$\\

It remains to study the case $\nu >0.$ The treatment of this case
is more involved and will be divided into several Lemmas. Consider
the setting of Proposition \ref{wx}, with $\nu
>0$. Since $M$ is compact, $\psi$ satisfies the equation
\begin{equation}\label{confkill}
\frac{2}{3}\Delta \psi  -q(R)\psi  +\frac{4(n-1)}{3(4n-1)}d X  =0.
\end{equation}
In the following Lemma we show that the
$S^{2}H\otimes\Lambda^{2}_{0}E$-component of $\psi$ is trivial and
that its $S^{2}H$-component is a solution of the twistor equation.
Recall that the map which associates to a Killing vector field $X$
on $(M, g)$  the $2$-form $\frac{2}{3\nu}(\nabla X)^{S^{2}H}$
(sometimes referred as the Hamiltonian form of $X$) is an
isomorphism (in particular, it is injective) from the space of
Killing vector fields to the space of solutions of the twistor
equation (see \cite{amp1}, Proposition 5.6). Moreover, the
Hamiltonian form of any Killing vector field on $(M, g)$  is an
eigenform of the Laplace operator $\Delta$, with eigenvalue $2\nu
(n+2)$ (see \cite{amp}, Theorem 2.7).

\begin{lem}\label{inc} The conformal Killing $2$-form $\psi $ is a section of
the direct sum bundle $S^{2}H\oplus S^{2}E$ and
\begin{equation}\label{s2h}
\psi ^{S^{2}H} =-\frac{2}{\nu (4n-1)}(\nabla X)^{S^{2}H}.
\end{equation}
In particular, $\psi^{S^{2}H}\neq 0$ unless $\psi$ is a Killing
$2$-form.
\end{lem}

\begin{proof}
Projecting the conformal-Killing equation
(\ref{confkill}) on $S^{2}H$ and using (\ref{constante}) we get
\begin{equation}\label{23}
\frac{2}{3}\Delta (\psi^{S^{2}H} )-4\nu \psi^{S^{2}H}
+\frac{8(n-1)}{3(4n-1)}( \nabla X )^{S^{2}H}=0.
\end{equation}
Define an operator
$$
T:\Gamma (Q)\rightarrow\Gamma (Q), \quad T(u):= \frac{2}{3}\Delta
u -4\nu u.
$$
The operator $T$ obviously preserves the eigenbundle $E_{2\nu
(n+2)}(Q)$ of $\Delta :\Gamma (Q)\rightarrow \Gamma (Q)$,
corresponding to the eigenvalue $2\nu (n+2)$, as well as its
orthogonal complement $E_{2\nu (n+2)}(Q)^{\perp}$, taken with
respect to the Euclidian metric of $\Gamma (Q)$ defined by the
metric $\langle \cdot ,\cdot \rangle$ of $Q$, followed by
integration over $M$. If we write $\psi^{S^{2}H}
=\psi_{1}+\psi_{2}$, with $\psi_{1}\in E_{2\nu (n+2)}(Q)$ and
$\psi_{2}\in E_{2\nu (n+2)}(Q)^{\perp}$ then, clearly,
$T(\psi_{2}) =0$ (because $(\nabla X)^{S^{2}H}\in E_{2\nu
(n+2)}(Q)$, as mentioned above). On the other hand, the
eigenvalues of $\Delta$ on $\Gamma (Q)$ are bigger or equal to
$2\nu (n+2)$ (see \cite{semqk}, Proposition 4.4). Since $6\nu <
2\nu (n+2)$, $6\nu$ cannot be an eigenvalue of $\Delta$ on $\Gamma
( Q )$. We deduce that $\psi_{2}=0$ and  $\psi^{S^{2}H}\in E_{2\nu
(n+2)}(Q).$ Using (\ref{23}), we now easily get (\ref{s2h}). When
$\psi$ is conformal-Killing but not Killing, $X$ is non-trivial
and $\psi^{S^{2}H}\neq 0$ (see the comments above).

In a similar way, we prove that $\psi$ is a section of $S^{2}H
\oplus S^{2}E$, i.e. that $\psi^{S^{2}H\otimes
\Lambda^{2}_{0}E}=0.$ For this, notice that, since $X$ is Killing,
$\nabla X\in \Gamma (S^{2}H\oplus S^{2}E)$ (see \cite{kn}, page
247). Projecting (\ref{confkill}) on
$S^{2}H\otimes\Lambda^{2}_{0}E$ and using (\ref{constante}) we get
\begin{equation}\label{proj}
\frac{2}{3}\Delta (\psi^{S^{2}H\otimes \Lambda_{0}^{2}E}) -2\nu (n+2)\psi^{S^{2}H\otimes\Lambda_{0}^{2}E}=0.
\end{equation}
The eigenvalues of $\Delta$ on $\Gamma (S^{2}H\otimes
\Lambda^{2}_{0}E)$ are greater or equal to $4\nu (n+1)$ (see
\cite{semqk}, Proposition 4.4). Using (\ref{proj}), we deduce that
$\psi^{S^{2}H\otimes \Lambda^{2}_{0}E}=0$ when $n>2$ (because in
this case $4\nu (n+1)> 3\nu (n+2)$). It remains to see what
happens when $n=2.$ When $n=2$ $(M, g)$ is a symmetric space (see
\cite{lebrun1}, Theorem 5.4) and, as mentioned in Section
\ref{prel}, $\psi$ can be written as a sum $a_{1}\psi_{1}+\cdots +
a_{k}\psi_{k}$ where $\psi_{i}$ are conformal-Killing $2$-forms
and also eigenforms of Laplacian $\Delta$, with eigenvalues, say,
$\nu_{i}.$ Some of the $\psi_{i}$'s are Killing $2$-forms, hence
parallel \cite{qK}. Others are conformal-Killing, but not Killing
$2$-forms, and for them the corresponding eigenvalues $\nu_{i}$
are all equal to $8 \nu$, from (\ref{s2h}). In particular, $\psi$
is the sum of a parallel $2$-form and an eigenform of $\Delta$
with eigenvalue $8\nu .$ Using (\ref{proj}) with $n=2$ we deduce
that $\psi^{S^{2}H\otimes \Lambda^{2}_{0}E}=0$.
\end{proof}

In order to prove Proposition \ref{wx}, it will be useful to write
the conformal-Killing equation on $2$-forms in an alternative way.
We do this in the next Lemma, whose proof is similar to the proof
of Proposition $19$ of \cite{cald}.

\begin{lem}\label{ecd}
The conformal-Killing $2$-form $\psi$ satisfies
\begin{equation}\label{cheie}
\nabla_{Y}\psi =\frac{1}{4n-1}\left( X\land
Y+\sum_{k=1}^{3}J_{k}X\land J_{k}Y -\sum_{k=1}^{3}\omega_{k}(X,
Y)\omega_{k}\right)
\end{equation}
for any $Y\in TM$. Here $\{ J_{1}, J_{2}, J_{3}\}$ is an
admissible basis of $Q$ with associated K\"{a}hler forms
$\omega_{1}, \omega_{2}, \omega_{3}$, and, as before, $X =\delta
(\psi )$ is the codifferential of $\psi$.
\end{lem}

\begin{proof}
In order to verify (\ref{cheie}), we will show that
\begin{equation}\label{dpsi}
d\psi = -\frac{3}{4n-1}\left( J_{1}X\land\omega_{1}+J_{2}X\land\omega_{2}+J_{3}X\land \omega_{3}\right) .
\end{equation}
Indeed, once we have (\ref{dpsi}), we replace it
into the conformal-Killing equation (\ref{confkil}) and we get
(\ref{cheie}).
To prove (\ref{dpsi}), we define a $3$-form
$$
\beta : = d\psi +\frac{3}{4n-1}
\left( J_{1}X\land\omega_{1}+J_{2}X\land\omega_{2}+J_{3}X\land \omega_{3}\right)
$$
and we show that
\begin{equation}\label{beta}
\beta (Y, JV, U) +  \beta (Y, V, JU) = 0,\quad \forall Y, U, V\in
TM,
\end{equation}
for a $J\in Z$, which, without loss of generality, can be taken
to be $J_{1}.$ (It is easy to check that a $3$-form
with the symmetries (\ref{beta}) must be zero; this implies
(\ref{dpsi}) and our claim).

To show (\ref{beta}), we evaluate $d\psi (Y, JV,U)+d\psi (Y, V,
JU)$ using the conformal-Killing equation, written in the form
\begin{equation}\label{cke}
\frac{1}{3}i_{Y}d\psi = \frac{1}{4n-1}Y\land X+\nabla_{Y}\psi
,\quad\forall Y\in TM.
\end{equation}
Notice that
\begin{equation}\label{david}
(Y\land X)(JV, U)+(Y\land X)(V, JU) = - (\omega \land J X )(Y, JV,
U)-(\omega\land JX)(Y, V, JU),
\end{equation}
where $\omega :=\omega_{1}$ denotes the K\"{a}hler  form associated to $J$.
Also, using Lemma \ref{inc}, we can write
\begin{align*}
 (\nabla_{Y}\psi )(JV,U) +
(\nabla_{Y}\psi )(V, JU) &
=-\frac{2}{\nu (4n-1)}\nabla_{Y} ( \nabla X)^{S^{2}H}(JV, U)\\
&-\frac{2}{\nu (4n-1)}\nabla_{Y}(\nabla X)^{S^{2}H}(V,JU)\\
&+\nabla_{Y}( \psi^{S^{2}E})(JV, U)+\nabla_{Y}( \psi^{S^{2}E})(V, JU).\\
\end{align*}
But $\nabla_{Y}( \psi^{S^{2}E})$ is $J$-invariant, hence
orthogonal to $JV\land U+V\land JU$ (which is $J$-anti-invariant)
and therefore
$$
\nabla_{Y}( \psi^{S^{2}E})(JV, U)+\nabla_{Y}( \psi^{S^{2}E})(V,
JU) =0.
$$
Moreover, from the Konstant formula
and the expression (\ref{rg}) of $R^{g}$
we deduce that
$$
\nabla_{Y}(\nabla X)^{S^{2}H}= R^{g}(Y, X)^{S^{2}H} = -
\frac{\nu}{2}\sum_{i=1}^{3}\omega_{i}(Y, X)\omega_{i}.
$$
Combining the above relations, we obtain
\begin{align*}
 (\nabla_{Y}\psi )(JV,U) + (\nabla_{Y}\psi )(V, JU)
&=\frac{2}{4n-1}\left( \omega_{2}(V, U)\omega_{3}(Y, X)
-\omega_{3}(V, U)\omega_{2}(Y, X)\right)\\
&= -\frac{1}{4n-1}\left( \alpha (Y, JV, U)+\alpha (Y, V, JU),
\right)
\end{align*}
where $\alpha\in \Omega^{3}(M)$ is defined by
$$
\alpha := J_{2}X\land \omega_{2}+J_{3}X\land\omega_{3}.
$$
Combining the above equality with (\ref{cke}) and (\ref{david}),
we obtain (\ref{beta}) and our claim.
\end{proof}

Using Lemma \ref{inc}, we can write our conformal-Killing $2$-form
$\psi$ as
$$
\psi=-\frac{2}{\nu (4n-1)}(\nabla X)^{S^{2}H}+ \frac{4}{\nu (4n-1)}(\nabla X)^{S^{2}E} +u,
$$
where $u$ is a section of $S^{2}E$. Equation (\ref{cheie})
written in terms of $u$
becomes
\begin{equation}\label{ad}
\nabla_{Y}u =-\frac{4}{\nu (4n-1)}W(Y, X)\quad\forall  Y\in TM .
\end{equation}

We conclude the proof of Proposition \ref{wx} with the following Lemma.

\begin{lem} Let $(M, g)$ be a compact quaternionic-K\"{a}hler manifold,
with positive reduced scalar curvature $\nu >0$ and
quaternionic-Weyl tensor $W$. Let $X$ be an arbitrary vector field
on $M$. Any section $u\in \Gamma (S^{2}E)$ which satisfies
(\ref{ad}) is parallel. In particular, $W(X,\cdot )=0.$
\end{lem}

\begin{proof}
Consider the exterior derivative $du$, written in the form
\begin{equation}
(du)(Z_{0}, Z_{1}, Z_{2})= ( \nabla_{Z_{0}} u)(Z_{1}, Z_{2}) -(\nabla_{Z_{1}}u)(Z_{0},
Z_{2}) +( \nabla_{Z_{2}}u)(Z_{0}, Z_{1}),
\end{equation}
where $Z_{0}, Z_{1}, Z_{2}\in TM.$
Using (\ref{ad}), we get
\begin{align*}
(du)(Z_{0}, Z_{1}, Z_{2}) &= \frac{4}{\nu (4n-1)}\left( - W(Z_{0},
X , Z_{1}, Z_{2}) +
W(Z_{1}, X, Z_{0}, Z_{2})\right)\\
& -\frac{4}{\nu (4n-1)}W(Z_{2}, X, Z_{0}, Z_{1}).
\end{align*}
The symmetries of the curvature tensor
$W$ imply that $du =0.$
Similarly, we can write the codifferential of $u$ in the form
$$
\delta u = -\sum_{i=1}^{4n}(\nabla_{e_{i}}u)(e_{i},\cdot )
=\frac{4}{\nu (4n-1)}W(e_{i}, X)(e_{i}) = 0,$$ because $W$ is in
the kernel of the Ricci contraction. We have proved that $u$ is
harmonic. Recall that the second Betti number $b_{2}(M)$ of $(M,
g)$ is zero, unless $(M,g)$ is isomorphic to the Grassmannian
$\mathrm{Gr}_{2}({\mathbb C}^{n+2})$ of complex $2$-planes in
$\mathbb{C}^{n+2}$, with its standard quaternionic-K\"{a}hler
metric \cite{lebrun1}; moreover, the space of harmonic $2$-forms
on $\mathrm{Gr}_{2}({\mathbb C}^{n+2})$ is one dimensional,
generated by the K\"{a}hler form, which is a parallel section of
$S^{2}E$. This proves that $u$ is actually parallel. Our claim
follows.
\end{proof}

\section{Killing vector fields and the quaternionic-Weyl tensor}\label{kerw}

In this Section we conclude the proof of Theorem \ref{main}.
We do this by proving Proposition \ref{kil}
stated below, which in turn relies on the following Theorem of Obata
(see \cite{obata}, Theorem $C$).

\begin{thm}\label{obatat} Let $(N^{2n}, J,g)$ be a complete,
connected and simply
connected K\"{a}hler manifold.
Suppose there is a non-constant smooth function $f$ on $N$ which satisfies
the Obata's equation
\begin{align*}
4\nabla^{2}(df)(Y, U, V)&= -2df(Y)g(U, V) - df(U)g(Y, V) - df(V) g(Y, U)\\
& + df(JU) \omega (Y, V)+ df(JV) \omega (Y, U),
\end{align*}
for any vector fields $Y, U, V\in {\mathcal X}(N)$, where $\nabla$
is the Levi-Civita connection and $\omega$ is the K\"{a}hler form.
Then $(N, J, g)$ is isometric to $(\mathbb{C}P^{n},
g_{\mathrm{FS}})$, where $g_{\mathrm{FS}}$ is the Fubini-Study
metric of constant holomorphic sectional curvature equal to one.
\end{thm}

\begin{rem}\label{rem1}{\rm  It is easy to verify
that the Hamiltonian function of any Killing
vector field on $(\mathbb{C}P^{n}, g_{\mathrm{FS}})$
satisfies the Obata's equation.
Conversely,
Theorem \ref{obatat} implies that the existence
of a {\it single} Killing vector field on a
complete, connected and simply connected
K\"{a}hler manifold,
whose Hamiltonian function is a solution of
the Obata's equation, insures that the K\"{a}hler manifold is
isometric to $(\mathbb{C}P^{n}, g_{\mathrm{FS}})$.}
\end{rem}

Proposition \ref{kil} below concerns compact quaternionic-K\"{a}hler
manifolds with positive scalar curvature.
Without loss of generality, we will normalise the
quaternionic-K\"{a}hler metric to have reduced scalar curvature $\nu =1.$
We shall denote by $g_{\mathrm{can}}:=g_{\mathrm{can}}(1)$ the
standard quaternionic-K\"{a}hler metric of $\mathbb{H}P^{n}$,
normalized in this way. The main result of this Section is the
following.

\begin{prop}\label{kil} Let $(M, g)$ be a
compact, quaternionic-K\"{a}hler manifold, of dimension $4n\geq 8$
and reduced scalar curvature $\nu =1.$ Suppose there is a
non-trivial Killing vector field $X$ on $M$ such that $W(X, \cdot
)=0$, where $W$ is the quaternionic-Weyl tensor. Then $W=0$ and
$(M,g)$ is isometric to $(\mathbb{H}P^{n}, g_{\mathrm{can}})$.
\end{prop}

\begin{rem}{\rm The idea of the proof of Proposition
\ref{kil} is to
show that the Hamiltonian function $f^{X}$ of the natural
lift $X^{Z}$ of $X$ to the K\"{a}hler-Einstein
twistor space $(Z, \bar{g}, {\mathcal J})$
of $(M, g)$ satisfies the Obata's equation stated above.
Theorem \ref{obatat} implies
that $(Z, \bar{g}, {\mathcal J})$ is isomorphic
to $(\mathbb{C}P^{n}, g_{FS})$ and then $(M, g)$ is isomorphic to
$(\mathbb{H}P^{n}, g_{\mathrm{can}})$. Details are as follows.}
\end{rem}

Since $X$ is Killing, its natural lift $X^{Z}$ on
$(Z, \bar{g}, {\mathcal J})$ is a Killing, real holomorphic
vector field, and its value at a point $J\in Z$ is
\begin{equation}\label{lift}
X^{Z}_{J} = \bar{X}_{J}+[\nabla X, J].
\end{equation}
In (\ref{lift}) (and in the following considerations), the bar
over a tangent vector on $M$
denotes its horisontal lift to $Z$, using the
Levi-Civita connection $\nabla$ of $g$, acting
on the twistor bundle $\pi :Z\to M$. The comutator
$[\nabla X, J]$ is seen as a tangent
vertical vector of $Z$ at $J$, and can be alternatively written as
\begin{equation}\label{lift1}
[\nabla X, J] = - 2\mathcal J(\tilde{A})_{J},
\end{equation}
where $A:= (\nabla X)^{S^{2}H}$ and $\tilde{A}$ is the
induced vertical vector field on $Z$, defined, at a point
$J\in Z$, by
\begin{equation}\label{a}
\tilde{A}_{J}:= A-\langle A, J\rangle J,\quad\forall J\in Z.
\end{equation}
Therefore, the vector field $X^{Z}$ has the form
\begin{equation}\label{new}
X^{Z} =\bar{X} -2{\mathcal J}(\tilde{A}).
\end{equation}
Since $(Z, {\mathcal J}, \bar{g})$ is compact,
K\"{a}hler-Einstein, with positive scalar curvature
$k^{\prime}= 2(2n+1)(n+1)$, the vector field
$X^{Z}$ is Hamiltonian, i.e.
$$
X^{Z} = {\mathcal J}\mathrm{grad}_{\bar{g}}\left( f^{X}\right)
$$
where
$$
f^{X}: = -\frac{1}{2(n+1)}\mathrm{trace}_{\bar{g}}\left( {\mathcal
J} \bar{\nabla}X^{Z}\right)
$$
is the Hamiltonian function of $X^{Z}$ and $\bar{\nabla}$
denotes the Levi-Civita connection of $\bar{g}.$
The way $f^{X}$ is related to the Hamiltonian
$2$-form of $X$ is explained in \cite{amp}.

In order to prove that $f^{X}$ satisfies the Obata's equation, we
need to calculate the second covariant derivatives of $\bar{X}$
and $\tilde{A}$ with respect to $\bar{\nabla}$, see relation
(\ref{new}). There are two types of tangent vectors on $Z$:
horisontal (i.e. which belong to the horisontal bundle determined
by $\nabla$, seen as a connection on the twistor bundle $\pi
:Z\rightarrow M$) and vertical. From Remark \ref{rem1}, the
Obata's equation for $f^{X}$ is satisfied when all arguments are
vertical (the twistor lines being totally geodesic and isomorphic,
as K\"{a}hler manifolds, to $(\mathbb{C}P^{1}, g_{\mathrm{FS}})$).
For this reason we shall not include, in Lemma \ref{1} and Lemma
\ref{2} below, the expressions of the second covariant derivatives
of $\bar{X}$ and $\tilde{A}$ when all arguments are vertical.

We begin with the following Lemma on the Levi-Civita connection
$\bar{\nabla} .$ This is by no means new, but we state it to fix
notations and conventions which will be useful later on in the
proof of Lemma \ref{1} and Lemma \ref{2}. We shall denote by
$\bar{\omega} =\bar{g}({\mathcal J}\cdot ,\cdot )$ the K\"{a}hler
form on $Z.$ As usual, for an admissible basis $\{ J_{1}, J_{2},
J_{3}\}$ of $Q$, $\omega_{1}$, $\omega_{2}$ and $\omega_{3}$ will
denote the corresponding K\"{a}hler forms.

\begin{lem}\label{lc}
For any vector fields $Y, V\in {\mathcal X}(M)$ and sections
$B, C\in \Gamma (Q)$,
\begin{align*}
\bar{\nabla}_{\bar{Y}}\bar{V} &= \overline{\nabla_{Y}V} -
\frac{1}{2}\left( \omega_{2}(Y,V)J_{3}
-\omega_{3}(Y,V)J_{2}\right)\\
\bar{\nabla}_{\tilde{B}}\bar{V}& = -
\frac{1}{2}{\mathcal J}( \tilde{B}\cdot V)\\
\bar{\nabla}_{\bar{Y}}\tilde{C} &
= -\frac{1}{2}{\mathcal J}( \tilde{C}\cdot Y)
+\widetilde{\nabla_{Y}C}\\
\bar{\nabla}_{\tilde{B}}\tilde{C}& = -\langle C, J\rangle \tilde{B}.
\end{align*}

The above expressions are evaluated at a point $J\in Z_{p}$, $\{ J
=J_{1}, J_{2}, J_{3}\}$ is an admissible basis of $Q$,
$\tilde{B}$, $\tilde{C}$ are vertical vector fields on $Z$ defined
as in (\ref{a}), $\tilde{B}\cdot V$ (and similarly for
$\tilde{C}\cdot Y$) is an horisontal vector field on $Z$, which at
$J$ is the horisontal lift of $\tilde{B}_{J}(V_{p})\in T_{p}M$
(here and below a tangent vertical vector $\tilde{B}_{J}\in
T_{J}Z_{p}$ is viewed also as an endomorphism of $T_{p}M$ and
$\tilde{B}_{J}(V_{p})$ denotes its action on $V_{p}\in T_{p}M$).

\end{lem}

\begin{rem}{\rm
In the same framework, the various Lie brackets
of basic and vertical vector fields
on the twistor space of a conformal-Weyl self-dual $4$-manifold
(i.e. a self-dual $4$-manifold together
with a fixed Weyl connection) were calculated in \cite{pgaud},
Appendix A. The same
formulas hold true also in the quaternionic-K\"{a}hler context,
with the Weyl connection replaced by the Levi-Civita
connection of the quaternionic-K\"{a}hler metric.}
\end{rem}

We now calculate the second covariant derivatives of $\bar{X}$ and
$\tilde{A}$ as follows.

\begin{lem}\label{1}
The second covariant derivative $\bar{\nabla}^{2}(\bar{X})$ at a
point $J\in Z_{p}$ has the following expression: for any tangent
vectors $B, C\in T_{J}Z_{p}$, $Y, U,V\in T_{p}M$ and admissible
basis $\{ J =J_{1}, J_{2}, J_{3}\}$ of $Q$,
\begin{align*}
\bar{g}\left(\bar{\nabla}^{2}(\bar{X})(\bar{Y},\bar{U}),
\bar{V}\right) &= \frac{1}{4}\left( \omega_{2}(Y,
V)\omega_{2}(X,U)+
\omega_{3}(Y, V)\omega_{3}(X,U)\right)\\
&+\frac{1}{4} \left(\omega_{2}(X, V)\omega_{2}(Y, U)+ \omega_{3}(X, V)\omega_{3}(Y, U)\right)\\
&+\frac{1}{2}\left(\omega_{2}(X,Y)\omega_{2}(U,V)+\omega_{3}(X,Y)\omega_{3}(U,V)\right)\\
\end{align*}
\begin{align*}
&+g( (X\land Y )^{S^{2}E}(U), V) +\frac{1}{2}\bar{\omega}(
\bar{X}, \bar{Y})\bar{\omega}(\bar{U},\bar{V});\\
\bar{g}\left( \bar{\nabla}^{2}(\bar{X})(\bar{Y}, \bar{U}),
B\right)&
= \frac{1}{2}\left( g(B( \nabla_{Y}X), JU) +g(B( \nabla_{U}X), JY)\right)\\
\bar{g}\left( \bar{\nabla}^{2}(\bar{X})(\bar{Y},  B),
\bar{U}\right)&
= -  \langle A, J \rangle g(B(Y), U) +\bar{g}(\tilde{A}, B) \bar{\omega}(\bar{Y}, \bar{U})\\
\bar{g}\left(\bar{\nabla}^{2}(\bar{X})(B, \bar{Y}),
\bar{U}\right)& =
\bar{g}(\tilde{A}, B) \bar{\omega}( \bar{Y}, \bar{U}) -\langle A, J\rangle g(B(Y), U)\\
\bar{g}\left(\bar{\nabla}^{2}(\bar{X})(B, \bar{Y}), C\right)
& = -\frac{1}{4}\left( \bar{\omega} (\bar{X},\bar{Y})\bar{\omega}(B, C)
+\bar{g}(\bar{X}, \bar{Y}) \bar{g}(B, C)\right) \\
\bar{g}\left(\bar{\nabla}^{2}(\bar{X})({B}, {C}), \bar{Y}\right)&
= -\frac{1}{4}\left(\bar{g}(\bar{X}, \bar{Y})\bar{g}({B}, {C})
 + \bar{\omega}(\bar{X},\bar{Y})\bar{\omega}({B}, {C})\right) \\
\bar{g}\left(\bar{\nabla}^{2}(\bar{X})(\bar{Y},B), C\right)&
= -\frac{1}{2}\bar{g}(B, C)\bar{g}(\bar{X}, \bar{Y}).\\
\end{align*}
\end{lem}

\begin{proof}
The proof uses Lemma \ref{lc} and is a straightforward calculation.
The condition $W(X, \cdot)=0$ comes into the picture
by means of the Konstant formula
$$
\nabla_{Y}(\nabla X) =R^{g}(Y, X) =
\frac{1}{4}\left( X\land Y +\sum_{i=1}^{3}J_{i}X\land J_{i}Y\right)
+\frac{1}{2}\sum_{i=1}^{3}\omega_{i}(X,Y)J_{i}.
$$
\end{proof}

In a similar way, we prove the following Lemma.

\begin{lem}\label{2} With the notations of Lemma \ref{1},
the second covariant derivative $\bar{\nabla}^{2}(\tilde{A})$
has the following expression at $J$:
\begin{align*}
\bar{g}\left( \nabla^{2}(\tilde{A})(\bar{Y},
\bar{U}),\bar{V}\right)
& = \frac{1}{4}\left( \omega_{2}(Y, X)\omega_{3}(U, V) -\omega_{2}(U, V)\omega_{3}(Y, X)\right)\\
& +\frac{1}{4}\left( \omega_{2}(U, X)\omega_{3}(Y, V)-\omega_{2}(Y, V)\omega_{3}(U, X)\right)\\
\bar{g}\left( \nabla^{2}(\tilde{A})(\bar{Y}, \bar{U}),B\right)&
= -\frac{1}{2}\left( g(B(U), \nabla_{Y}X)+\langle A, J\rangle g(B(Y), JU)\right)\\
&-\frac{1}{4}\left( \bar{\omega}(\bar{Y}, \bar{U}) \bar{\omega}(B,
\tilde{A})+\bar{g}(\bar{Y}, \bar{U})\bar{g}(\tilde{A},
B)\right)\\
\end{align*}
\begin{align*}
\bar{g}\left( \nabla^{2}(\tilde{A})(\bar{Y}, B),\bar{U}\right)&=
\frac{1}{4}\left( g(A(U), B(Y)) -\langle A, J\rangle g(B(Y),
JU)\right)\\
&-\frac{1}{2}\langle A, J\rangle g(B(Y), JU)\\
\bar{g}\left( \nabla^{2}(\tilde{A})(\bar{Y}, B),C\right)& =
\frac{1}{4}\left(\bar{g}(\bar{X}, \bar{Y})\bar{\omega}(B, C)
-\bar{\omega}(\bar{X}, \bar{Y})\bar{g}(B, C)\right)\\
\bar{g}\left( \nabla^{2}(\tilde{A})(B, \bar{Y}),\bar{U}\right)&
= -\frac{1}{2}\langle A, J \rangle g(B(Y), JU)\\
\bar{g}\left( \nabla^{2}(\tilde{A})(B, \bar{Y}),C\right)&
= \frac{1}{4}\left( \bar{g}(\bar{X}, \bar{Y})\bar{\omega}(B, C)
- \bar{\omega}(\bar{X}, \bar{Y})\bar{g}(B, C)\right)\\
\bar{g}\left( \bar{\nabla}^{2}(\tilde{A})({B}, {C}), \bar{Y}\right)&
= 0.
\end{align*}
\end{lem}

The following Lemma concludes
the proof of Proposition \ref{kil}.

\begin{lem}\label{ordine}
Consider the setting of Proposition \ref{kil}.
Then the Hamiltonian function
$f^{X}$ of the natural lift $X^{Z}$ of the Killing
vector field $X$ to the
twistor space $(Z,\bar{g}, {\mathcal J})$ satisfies the
Obata's differential equation.
In particular, $(Z, \bar{g}, {\mathcal J})$ is isomorphic to
$(\mathbb{C}P^{2n+1},g_{\mathrm{FS}})$ and $(M, g)$ is isomorphic
to $(\mathbb{H}P^{n}, g_{\mathrm{can}})$.
\end{lem}

\begin{proof}
It is straightforward (from Lemma \ref{1} and Lemma \ref{2}) to
check that $f^{X}$ satisfies the Obata's equation. Probably the
most involved computation is to check the Obata's equation for
$f^{X}$ when the first two arguments of
$\bar{\nabla}^{2}\left(df^{X}\right)$ are horisontal and the third
is vertical. To do this calculation, we write (with the notations
of the previous Lemmas) at $J$:
\begin{align*}
\bar{\nabla}^{2} (df^{X})(\bar{Y}, \bar{U}, B) & = \bar{g}\left(
\bar{\nabla}^{2}(\bar{X})(\bar{Y}, \bar{U}), JB\right)
-2\bar{g}\left( \bar{\nabla}^{2}(\tilde{A})(\bar{Y},\bar{U}), B\right)\\
&= \frac{1}{2}\langle B, J_{2}\rangle\left( (\nabla X)(Y, J_{2}U)
- (\nabla X)(U, J_{2}Y)\right)\\
&+\frac{1}{2} \langle B, J_{3}\rangle \left( (\nabla X)(Y, J_{3}U)
-(\nabla X)(U,J_{3}Y)\right)\\
& +\langle A, J\rangle \left( \omega_{2}(Y, U)\langle B,
J_{3}\rangle -\omega_{3}(Y, U)\langle B, J_{2}\rangle
\right)\\
& +\frac{1}{2}\left(\bar{\omega}(\bar{Y}, \bar{U}) \bar{\omega}(B,
\tilde{A}) +\bar{g}(\bar{Y},\bar{U}) \bar{g}(\tilde{A}, B)\right)
.
\end{align*}
On the other hand, since $Y\land J_{2}U -U\land J_{2}Y$ and
$Y\land J_{3}U - U\land J_{3}Y$ is $J_{2}$ (respectively, $J_{3}$)
anti-invariant,
\begin{align*}
( \nabla X)(Y, J_{2}U) -(\nabla X)(U, J_{2}Y) &= 2\left( \langle
A, J \rangle \omega_{3}(Y,U) +
\langle A, J_{3}\rangle \bar{\omega}(\bar{U}, \bar{Y})\right)\\
(\nabla X)(Y, J_{3}U) - ( \nabla X)(U, J_{3}Y)
 &= 2\left( \langle A, J \rangle \omega_{2}(U,Y) +
\langle A, J_{2}\rangle \bar{\omega}(\bar{Y},\bar{U})\right) .
\end{align*}
We now easily get
$$
\bar{\nabla}^{2}(df^{X})(\bar{Y}, \bar{U}, B) =\frac{1}{4} \left(
df^{X}({\mathcal J}B)\bar{\omega}(\bar{Y}, \bar{U}) - df^{X}(B)
\bar{g}(\bar{Y}, \bar{U})\right) ,
$$
i.e. $f^{X}$ satisfies the Obata's equation, when the first two
arguments are horisontal and the third is vertical. To conclude the
proof, it is enough to notice that $Z$ is simply connected (see
\cite{inventiones}, Theorem 6.6) and to apply Theorem \ref{obatat}
to $(Z, \bar{g}, {\mathcal J})$ and the function $f^{X}.$

\end{proof}

The proof of Theorem \ref{main} is thus completed.\\

It is worth to relate Proposition \ref{kil} to the theory on
gradient quaternionic vector fields, developed in \cite{gradient}
and \cite{gradient1}. Recall that a vector field on a quaternionic
manifold is quaternionic, if its flow preserves the quaternionic
bundle. On a complete quaternionic-K\"{a}hler manifold with
non-zero scalar curvature, any complete quaternionic vector field
is the sum of a Killing vector field and of a (complete) gradient
quaternionic vector field (see Proposition 4 of \cite{gradient1}).

\begin{cor} Let $(M, g)$ be a compact,
quaternionic-K\"{a}hler manifold of dimension $4n\geq 8$, with
quaternionic-Weyl tensor $W$. Suppose there is a non-trivial
quaternionic vector field $X$ of $(M, g)$, such that $W(X,\cdot
)=0.$ Then either $(M, g)$ is isometric to the standard
quaternionic projective space or it is Ricci-flat and its
universal cover has an Euclidian factor $(\mathbb{H}^{k}, g_{0})$
in its De-Rham decomposition.
\end{cor}

\begin{proof}
Suppose first that $X$ is Killing. We distinguish two cases: $\nu
>0$ and $\nu = 0$ (as already mentioned in the proof of Proposition \ref{wx},
there are no non-trivial Killing vector fields on $(M ,g)$ when
$\nu <0$). If $\nu >0$ then $(M, g)$ is isometric to the standard
quaternionic projective space (see Proposition \ref{kil}). If $\nu
= 0$ then $X$ is parallel (see \cite{der}, Theorem 1.84) and we
deduce that the universal cover of $(M, g)$ has a flat factor in
its De-Rham decomposition.

Next, suppose that $X$ is quaternionic, but not Killing. Again we
distinguish two cases: $\nu =0$ and $\nu \neq 0.$ If $\nu =0$ then
the universal cover of $(M,g)$ is complete, simply connected and
the claim has been proved in Theorem 1 of \cite{gradient1}. If
$\nu \neq 0$ then there is a gradient (non-trivial) quaternionic
vector field $\mathrm{grad}_{g}(f)$ on $(M, g)$. The existence of
such a vector field on $(M, g)$ implies again that $(M, g)$ is
isometric to the standard quaternionic projective space (see
\cite{gradient}, \cite{gradient1} and also \cite{lebrunijm}). In
fact, the proof of this statement when the scalar curvature of
$(M, g)$ is positive (see \cite{gradient}, Theorem 1) consists in
showing that the pull-back of $f$ to the total space $(S, g_{S})$
of the $3$-Sasaki bundle of $(M, g)$ (with $g$ normalised such
that $\nu =1$) satisfies an equation considered in Theorem A of
\cite{obata} and then to deduce, using the theory developed in
\cite{obata}, that $(S, g_{S})$ is isometric to the Euclidian
sphere of radius two. These considerations imply that $(M, g)$ is
isometric to $(\mathbb{H}P^{n}, g_{\mathrm{can}})$, see
\cite{gradient}. Alternatively, one could have noticed  that the
pulled back $\bar{f}:= \pi^{*}f$ of $f$ to the twistor space
$(Z,\bar{g}, {\mathcal J})$ of $(M, g)$ is the Hamiltonian
function of a Killing vector field on $(Z, \bar{g}, {\mathcal J})$
(see \cite{amp}, Proposition 3.1) and check instead that $\bar{f}$
satisfies the Obata's equation of Theorem \ref{obatat}.
\end{proof}

For Wolf spaces, the proof of  Proposition \ref{kil} can be considerably
simplified, by proving the following Lemma.

\begin{lem}
Let $(M, g)$ be a non Ricci-flat quaternionic-K\"{a}hler manifold,
with quaternionic Weyl tensor $W$. Suppose that $W(X,\cdot )=0$
for a non-trivial (not necessarily Killing) vector field. Then the
holonomy algebra of $(M, g)$ is the entire $\mathrm{sp}(1)\oplus
\mathrm{sp}(n).$ In particular, if $(M, g)$ is a Wolf space, then
it is necessarily isometric to the quaternionic projective space,
with its standard quaternionic-K\"{a}hler structure.
\end{lem}

\begin{proof} Since $(M, g)$ is non Ricci-flat,
the holonomy algebra $\mathrm{hol}(g)$ of $(M, g)$ contains the
$\mathrm{sp}(1)$-factor of $\mathrm{sp}(n)\oplus \mathrm{sp}(1).$
In order to prove our claim, we need to show that also $(Y\land
U)^{S^{2}E}$ belongs to $\mathrm{hol}(g)$, for any pair of tangent
vectors $Y$ and $U$. Let $\{ J_{1}, J_{2}, J_{3}\}$ be an
admissible basis of $Q$. Since $W(X,\cdot )=0$, also $W(J_{i}X,
\cdot )=0$, because $W$, viewed as a vector-valued $2$-form, is
$Q$-hermitian. Recall that the value of the curvature $R^{g}$ on
any pair of tangent vectors belongs to the holonomy algebra. Using
$\mathrm{sp}(1)\subset \mathrm{hol}(g)$, we deduce that
\begin{equation}\label{a1}
R^{g}(X, V)^{S^{2}E}= -\nu (X\land V)^{S^{2}E}\in\mathrm{hol}(g)
\end{equation}
and
\begin{equation}\label{a2}
R^{g}(J_{i}X, V)^{S^{2}E}= -\nu (J_{i}X\land V)^{S^{2}E}\in\mathrm{hol}(g)
\end{equation}
for any tangent vector $V\in TM.$ It follows that if $Y$ or $U$
belong to the vector space ${\mathcal V}:= \mathrm{Span}\{ X,
J_{1}X, J_{2}X, J_{3}X\}$, then $(Y\land U)^{S^{2}E}$ belongs to
$\mathrm{hol}(g).$ It remains to show that $(Y\land U)^{S^{2}E}$
belongs to the holonomy algebra when both $Y$ and $U$ are
orthogononal to ${\mathcal V} .$ Take such two tangent vectors $Y$
and $U$. Notice that, since both $(X\land Y)^{S^{2}E}$ and
$(X\land U)^{S^{2}E}$ belong to $\mathrm{hol}(g)$, also their Lie
bracket, which is equal to
$$
[(X\land Y)^{S^{2}E}, (X\land U)^{S^{2}E}] =
\frac{1}{16}\sum_{i,j=1}^{3} g(J_{i}Y, J_{j}U)J_{i}X\land J_{j}X
+\frac{1}{4} g(X, X) (Y\land U)^{S^{2}E},
$$
belongs to $\mathrm{hol}(g)$, as well as
the $S^{2}E$-part of this Lie bracket. Using (\ref{a1}) and (\ref{a2})
we get our claim.

\end{proof}

\section{The dimension of the space of conformal-Killing
$2$-forms}

It is known that the dimension $\mathrm{ck}_p(M)$ of the space of
conformal-Killing $p$-forms on a Riemannian manifold $(M, g)$ is
always finite whether $M$ is compact or not \cite{semc}. In this
Section we determine $\mathrm{ck}_{2}(M)$, when $(M,g)$ is
quaternionic-K\"{a}hler and compact. We begin with the following
considerations on quaternionic-K\"{a}hler manifolds with zero
scalar curvature.

\begin{rem}\label{dimensiunea}{\rm
Let $(M,g)$ be a compact quaternionic-K\"ahler manifold of zero
scalar curvature. Being Ricci-flat, $(M, g)$ has a finite
Riemannian covering $(T^{4q}\times \bar{M}, g_{0}\times {g}_{1})$,
where $T^{4q}$ is a $4q$-dimensional torus, with flat metric
$g_{0}$, and $(\bar{M}, {g}_{1})$ is compact and simply connected
(this is a result of Cheeger and Gromoll, see \cite{der}, page
169). Since $(M, g)$ is quaternionic-K\"{a}hler, $(\bar{M},
{g}_{1})$ is hyper-K\"{a}hler and can be decomposed into a
Riemannian product
$$
\bar{M} = S_{1}\times\cdots \times S_{l}
$$
where $S_{i}$ are hyper-K\"{a}hler, irreducible, of dimension
$4r_{i}$ and $\mathrm{Hol}(S_{i}) = \mathrm{Sp}(r_{i}).$ The Deck
group $G$ of the covering $T^{q}\times \bar{M}\rightarrow M$ is
included in the isometry group of $(T^{q}\times \bar{M},
g_{0}\times {g}_{1})$ and hence is a product group $G= H\times I$
(because the metric of $T^{4q}\times \bar{M}$ is a product
metric). Defining $F: = T^{4q}/ H$ we obtain a new Riemannian,
finite covering of $(M,g)$, isometric to
\begin{equation}\label{n}
N = F\times S_1\times\cdots  \times S_{l},
\end{equation}
with the following properties:
\begin{enumerate}

\item $F$ is a flat manifold finitely covered by a hyper-K\"ahler
torus;

\item the Deck group of the covering $N\rightarrow M$ is $I$. The isometric action
of $I$ on $N$ is the product action and is trivial on the first
factor $F$;

\item if $M$ has finite fundamental group, then
$N = S_{1}\times\cdots\times S_{l}$ (see \cite{der}, Corollary
6.67(a), page 168) and $M$ is hyper-K\"{a}hler if and only if
$M=N$, or, equivalently, if and only if $M$ is simply connected.

\end{enumerate}

Because each $S_i$ is simply connected $b_2(N)=b_2(F) +
\sum_{i=1}^l b_2 (S_i)$. Of course $b_2(F)$ equals the number of
parallel $2$-forms on $F$, while on each factor $S_i$ there are
exactly three parallel $2$-forms -- those coming from the
hyper-K\"ahler structure -- because otherwise the holonomy would
be strictly contained in $\mathrm{Sp}(r_i)$ (see also \cite{fuji},
Proposition 3.15). This implies that, for a compact
quaternionic-K\"ahler manifold with $\nu=0$ and {\it finite}
fundamental group, there are no parallel 2-forms in the subbundle
$(\Lambda_0^2E\otimes S^2H)\oplus S^2E$. In fact, since the
curvature of $\Lambda_0^2E\otimes S^2H$ vanishes, in the
terminology of \cite{semqk} we actually have
$b_{\textrm{expt,2}}(M)=0$ in this case.

A final observation before stating our next result  is that a
compact $4n$-dimensional hyper-K\"{a}hler manifold with holonomy
{\it equal} to $\mathrm{Sp}(n)$ is simply connected (see
\cite{der}, Lemma 14.21 and also \cite{fuji}, Remark 4.1). In
particular, a smooth finite quotient $\hat S_{i}=S_{i}/\Gamma$ of
a hyper-K\"{a}hler manifold $S_{i}$ from the decomposition
(\ref{n}) cannot be hyper-K\"ahler, unless $\Gamma$ is trivial. It
follows that the number of parallel $2$-forms on $\hat{S}_{i}$ is
at most one, and is one precisely when $\hat{S}_{i}$ is
K\"{a}hler. }
\end{rem}

\medskip

We now state the main result of this Section. We remark that our
treatment when the scalar curvature is negative is not complete;
however, all known (to us) {\it compact} examples of
quaternionic-K\"{a}hler manifolds with $\nu <0$  are locally
symmetric.

\begin{prop} Let $(M,g)$ be a compact quaternionic-K\"ahler
manifold of real dimension $4n\geq 8$ and reduced scalar curvature
$\nu$.

\begin{enumerate}

\item If $\nu>0$,

$$ \mathrm{ck}_2(M)= \left\{ \begin{array}{ll}
(n+1)(2n+3) & \textnormal{ if $M$ is standard $\mathbb H P^{n}$} \\
1           & \textnormal{ if $M$ is standard  $G_2(\mathbb C ^{n+2})$}\\
0           & \textnormal{ otherwise.}
\end{array} \right.
$$

\item If $\nu=0$,

 $$\mathrm{ck}_2(M)= b_2(F)+3l$$

if and only $\pi_1(M)=\pi_1(F)$ -- i.e. $M= N$ in Remark
\ref{dimensiunea}. Otherwise, $b_2(F)\leq \mathrm{ck}_2(M)\leq
b_2(F)+3l-2$ and examples can be constructed to show that every
possible value does occur.

\item If $\nu<0$ and the universal covering $\tilde M$ of $(M,g)$ is
symmetric,
$$ \mathrm{ck}_2(M)= \left\{ \begin{array}{ll}
1           & \textnormal{ if $\tilde M$ is the non-compact dual
of $\mathrm{Gr}_{2}({\mathbb C}^{n+2})$}\\
0           & \textnormal{ otherwise.}
\end{array} \right.
$$

\end{enumerate}

\end{prop}

\begin{proof}
Since any Killing $p$-form on a compact quaternionic-K\"ahler
manifold is parallel \cite{semqk}, we essentially have to count
parallel $2$-forms unless $M$ admits a conformal-Killing $2$-form
$\psi$ which is non-Killing. By Theorem \ref{main} this happens
only when $M$ is the standard $\mathbb{H}P^{n}$ in which case
$\psi$ is given as in Proposition 2 with $u=0$ (because
$b_2(\mathbb{H}P^{n})=0$) and the $S^2H$-component of $\psi$ is a
non-zero solution of the twistor equation. By Lemma 6.5 of
\cite{inventiones}, the composition
$$
\psi \rightarrow \psi^{S^{2}H}\rightarrow \delta (\psi^{S^{2}H})
$$
is an isomorphism from the space of conformal-Killing $2$-forms to
the space of Killing vector fields on $\mathbb{H}P^{n}$, so that
$\mathrm{ck}_2(\mathbb{H}P^{n})=\mathrm{dim}(\mathrm{Isom}(\mathbb{H}P^{n}))$.
In all other cases $\mathrm{ck}_2(M)$ is the dimension of the
space of parallel $2$-forms. We shall treat the cases $\nu >0$,
$\nu =0$ and $\nu <0$ separately, as follows.

Consider first the case $\nu >0.$ Recall that the space of
parallel $2$-forms on $M= \mathrm{Gr}_{2}({\mathbb C}^{n+2})$ is
generated by the K\"{a}hler form and is one dimensional. Moreover,
recall that $b_2(M)=0$ when $M$ is compact, with positive scalar
curvature and is not isometric to the Grassmannian
$\mathrm{Gr}_{2}({\mathbb C}^{n+2})$ (see \cite{lebrun1}). This
concludes the case $\nu >0.$

The case $\nu=0$ easily follows from Remark \ref{dimensiunea}.
Examples of Ricci-flat compact quaternionic-K\"{a}hler manifolds
$(M, g)$ with all possible values of $\mathrm{ck}_2(M)$ are
provided by products of finite quotients of $K3$-surfaces.

When $\nu <0$, certainly no solution of the twistor equation
exists on $(M,g)$ (see \cite{lili}, Theorem 9); the cohomology of
$M$ is the direct sum of the $\mathrm{sp}(1)$-invariant and
Exceptional cohomology (see \cite{semqk}, page 402); furthermore
$b_{\textrm{expt,2}}(M)=0$ at least for $n\geq 3$ (see
\cite{semqk}, Proposition 6.8). In any case, it follows from
relation (\ref{constante}) that any parallel $2$-form on $(M, g)$
is a section of $S^2 E$. Consider now the special case when the
universal covering $\tilde M=G^{*}/ K$ of $(M,g)$ is symmetric.
Parallel $2$-forms on $(M, g)$ lift to parallel $2$-forms on
$\tilde{M}$. These, in turn, are provided by  $2$-forms preserved
by the holonomy representation of $K$. In particular, parallel
$2$-forms on $\tilde{M}$ are in one to one correspondence with
parallel $2$-forms on the dual $G/ K$ of $\tilde{M}$. Our last
claim follows from the considerations we did in the case $\nu>0$.

\end{proof}

LIANA DAVID: Institute of Mathematics "Simion Stoilow" of the
Romanian Academy, Calea Grivitei nr. 21, Bucharest, Romania.\\
{\it E-mail address:} liana.david@imar.ro

\vspace{10pt}

MASSIMILIANO PONTECORVO: Dipartimento di Matematica, Universita' di
Roma Tre, L.go S.L. Murialdo 1, 00146, Roma, Italy.\\
{\it E-mail address:} max@mat.uniroma3.it

\end{document}